\documentclass{article}

\usepackage{amssymb}
\usepackage{url}
\newtheorem{theorem}{Theorem}
\newtheorem{lemma}{Lemma}

\begin{document}
\title{On Vaughan's approximation in restricted sets of arithmetic progressions}
\author{Claus Bauer}

\date{}

\maketitle
\begin{abstract}
We investigate the approximation to the number of primes in arithmetic progressions given by Vaughan (\cite{vau}). Instead of averaging the expected error term over all residue classes to modules in a given range, here we only consider subsets of arithmetic progressions that satisfy additional congruence conditions and provide asymptotic approximations.
 \end{abstract}

\section{Introduction}\setcounter{equation}{0}\setcounter{theorem}{0}\setcounter{lemma}{0}

The distribution of primes in arithmetic progressions is a long standing topic in analytic number theory. Defining
\begin{eqnarray*}
\vartheta(x,q,b)= \sum\limits_{p\leq x\atop p\equiv b(mod\,d)}\log\,p,
\end{eqnarray*}
the approximation $\vartheta(x,q,b)\sim\frac{x}{\phi(q)}$ is known to be true for small $q$ and $(a,q)=1.$ For larger $q,$
in \cite{barban}, \cite{hooley1}, and \cite{mont1} the quantity
\begin{eqnarray}\label{eq:barban}
&&\sum\limits_{d\leq Q}\sum\limits_{b=1\atop (b,d)=1}^{d}E^{2}(x,d,b),\\
&&E(x,d,b)= \vartheta(x,d,b)-\frac{x}{\phi(d)},\nonumber
\end{eqnarray} has been investigated for $Q\leq x.$
The most accurate asymptotic expression for (\ref{eq:barban}) was given in \cite{goldston} as follows:\\
\\
{\it For $Q<x,$ any $A>0,$ and a constant C,}
\begin{eqnarray}\label{eq:hoo}
\sum\limits_{d\leq Q}\sum\limits_{b=1\atop (b,d)=1}^{d}E^{2}(x,d,b)=
Qx\log\,Q+CQx+O\left(Q^{3/2}x^{1/2}+x^{2}(\log\,x)^{-A}\right).
\end{eqnarray}
In \cite{vau}, it was shown that for large $q$ the approximation $\frac{x}{\phi(q)}$ used in (\ref{eq:barban}) is not best possible. Setting
\begin{eqnarray*}
  &&  F_{R}(n)=\sum\limits_{r\leq R}\frac{\mu(r)}{\phi(r)}\sum\limits_{b=1\atop (b,r)=1}^{r}e\left(\frac{bn}{r}\right),\qquad
    \Delta(n)=\Lambda(n)-F_{R}(n),\\
    &&\rho(x,d,b)=\sum\limits_{n\leq x\atop n\equiv b(mod\,d)}F_{R}(n),
    \end{eqnarray*}the following approximation
was shown in \cite{vau}:
    \begin{theorem}\label{th:vau} For any positive integer $A$ and
$R=L^{A},$ there is
  \begin{eqnarray}\label{eq:vv}
\sum\limits_{d\leq Q}\sum\limits_{b=1}^{d}\left(\vartheta(x,d,b)-\rho(x,d,b)\right)^{2}
=Qx\log(x/R)-c_{0}Qx+O\left(QxR^{-1/2}+x^{2}(\log\,x)^{2}R^{-1}\right),
\end{eqnarray}
where
$
c_{0}=1+\gamma+\sum\limits_{p\geq 2}\frac{\log\,p}{p(p-1)}.
$
\end{theorem}
For $x/R<Q,$ the main term in (\ref{eq:vv}) is smaller than the main term in (\ref{eq:hoo}).\\
\\
We notice that in Theorem \ref{th:vau} the average is taken over all reduced residue classes $b$ mod $d.$ It is of interest to understand if one can get a better asymptotic when limiting the summation to reduced rest classes only. Theorem \ref{th:5} shows that this is indeed the case.
          \begin{theorem}\label{th:5} For a real number $x,$ a real number $R$ satisfying
$R=L^{G}$ for some $G>10$ and $Q=x(\log\,x)^{-B}$ for some $B>0,$ there is
  \begin{eqnarray*}
&&\sum\limits_{d\leq Q}\sum\limits_{b=1\atop (b,d)=1}^{d}\left(\vartheta(x,d,b)-\rho(x,d,b)\right)^{2}\nonumber\\
&&xQ\left(\log\,\left(x/R^{\left(2-\zeta^{-1}(2)\right)}\right)\right)
+xQ\left(-c_{1}+1
-\prod\limits_{p\geq 2}\left(1-\frac{1}{p(p-1)}\right)+\prod\limits_{p\geq 2}\left(1-\frac{1}{p^{2}}\right)c_{2}\right)\nonumber\\
&+&O\left(\frac{xQ}{R^{1/2}}+\frac{x^{2}(\log\,x)^{2}}{R}\right),
\end{eqnarray*} where
\begin{eqnarray*}
c_{1}&=&1+2\gamma+2\sum\limits_{p\geq 2}\frac{\log\,p}{p(p-1)}, \\
c_{2}&=&\gamma+\sum\limits_{p\geq 2}\frac{\log\,p}{p(p-1)}.
\end{eqnarray*}
    \end{theorem}
Comparing Theorems  \ref{th:vau} and \ref{th:5}, we see that the non-reduced residue classes with $b$ with $(r,b)>1$ indeed make a quantifiable contribution to the considered variance. This is different from the classical approach as in (\ref{eq:hoo}) where the contribution of the non-reduced rest-classed can be neglected. This observation is in alignment with the findings from \cite{fiorilli}. In \cite[Theorem 1.5]{fiorilli}, for a fixed $b,$ the average over $q$ of $\vartheta(x,d,b)-\rho(x,d,b)$ is considered. As in our case, it is shown that the non-reduced rest-classes make a non-negligible contribution.\\
\\
For the proof of Theorem \ref{th:5}, we will need the following Theorem that describes the distribution of the term $\Delta^{2}(n)$ in arithmetic progressions:
   \begin{theorem}\label{th:3}
For a real number $x,$ a non-negative integer $N,$ a positive integers $F$ and a squarefree integer $v,$ where $v\ll (\log\,x)^{F},$ and a real number $R\leq x^{1/3},$
\begin{eqnarray*}\label{eq:y1}
&&\sum\limits_{n\leq x\atop n\equiv N(mod\,v)}\Delta^{2}(n)
=\delta(N,v)\frac{x}{\phi(v)}\left(\log\,x-2\log\,R-c_{1}\right)+\frac{x}{v}\left(\log\,R+c_{2}\right)+\delta(N,v)x\frac{v}{\phi^{2}(v)}-\frac{x}{\phi(v)}\nonumber\\
&+&O\left(x\,exp\left(-cL^{1/2}\right)+\frac{x\tau(v)}{vR^{1/2}}+\frac{x}{\phi(v)R^{1/2}}+R^{2}\log\,R+\tau(v)R+\frac{x(\log\,v+\tau(v))}{v}\sum\limits_{r|v\atop r>R}1\right),
\end{eqnarray*} where $\delta(N,v)=1$ if either $N>0$ and $(N,v)=1$ or $N=0$ and $v=1,$ and $\delta(N,v)=0,$ otherwise, $c$ is a positive constant, and $c_{1}$ and $c_{2}$ are as defined in Theorem \ref{th:5}.
\end{theorem}

Further applying Theorem \ref{th:3}, we consider another variant of Theorem \ref{th:vau} where only average over such arithmetic progressions satisfying specific congruence conditions. We prove the following result:

\begin{theorem}\label{th:4} For a real number $x,$ a real number $R$ satisfying
$R=L^{G}$ for some $G>10$ and $Q=x(\log\,x)^{-B}$ for some $B>0,$ there is
  \begin{eqnarray*}
&&\sum\limits_{d\leq Q}\sum\limits_{b=1\atop (N-b,d)=1}^{d}\left(\vartheta(x,d,b)-\rho(x,d,b)\right)^{2}\nonumber\\
&=&xQ\prod\limits_{p\geq 2\atop (p,N)=1}\left(1-\frac{1}{p(p-1)}\right)\left(\log\,\left(x/R^{t(N)}\right)\right)\nonumber\\
&+&xQ\Biggl(-\prod\limits_{p\geq 2\atop (p,N)=1}\left(1-\frac{1}{p(p-1)}\right)c_{1}+
\prod\limits_{p\geq 2\atop (p,N)=1}\left(1-\frac{1}{(p-1)^{2}}\right)+\prod\limits_{p\geq 2}\left(1-\frac{1}{p^{2}}\right)c_{2}\nonumber\\
&&-\prod\limits_{p\geq 2}\left(1-\frac{1}{p(p-1)}\right)\Biggr)+O\left(\frac{xQ}{R^{1/2}}+\frac{x^{2}(\log\,x)^{2}}{R}\right),
\end{eqnarray*}
where $c_{1}$ and $c_{2}$ are as defined in Theorem \ref{th:5} and
\begin{eqnarray*}
t(N) &=&
\left(2-\prod\limits_{p\geq 2\atop (p,N)=1}\left(1-\frac{1}{p(p-1)}\right)^{-1}\prod\limits_{p\geq 2}\left(1-\frac{1}{p^{2}}\right)\right)\geq 1.
\end{eqnarray*}
    \end{theorem}
As the proofs of Theorem \ref{th:5} and \ref{th:4} are very similar, we will give the detailed proof of Theorems \ref{th:4}, but subsequently only provide a shortened version of the proof of Theorem \ref{th:5}.
We will use the abbreviation $L=\log\,x$ throughout this paper. $c$ denotes a positive absolute constant that can take different values at different occasions. $\epsilon$ denotes an arbitrarily small positive number and $\tau$ denotes the divisor function. We use the common abbreviations $e(x)=e^{2\pi ix}$ and
$\exp(x)=e^{x}.$
\section{Auxiliary Lemmas}\setcounter{equation}{0}\setcounter{theorem}{0}\setcounter{lemma}{0}
\begin{lemma}\label{le:th2vau}
There is a positive constant $c$ such that whenever $A$ is a fixed positive number, $R\leq (\log\,x)^{A},$ $\sqrt{x}\leq y\leq x,$ and $q\leq R,$ we have
\begin{eqnarray*}
\vartheta (y,q,a)-\rho(y,q,a)\ll y\,exp\left(-cL^{1/2}\right).
\end{eqnarray*}
\end{lemma}
{\it Proof:} \cite[Theorem 2]{vau}.
\begin{lemma}\label{le:hy} For $(a,q)=1,$ $q\leq x,\,q\leq R,$
\begin{eqnarray*}
\rho(x,q,a)\ll R +\frac{x\tau(q)}{q}.
\end{eqnarray*}
\end{lemma}
{\it Proof:} \cite[Theorem 1, (1.12)]{vau}.
\begin{lemma}\label{le:hilbert}
Let $v_{r}$ be any set of complex numbers and let $x_{r}$ be any set of real numbers distinct modulo $1.$
If $0<\delta=\min\limits_{r\ne s}|| x_{r}-x_{s}||,$ then for any real $t,$
\begin{eqnarray*}
\left| \sum\limits_{r\ne s}\bar{v_{r}}v_{s}\frac{\sin t(x_{r}-x_{s})}{\sin \,\pi(x_{r}-x_{s})}\right|\leq \sum\limits_{r}|v_{r}|^{2}\delta^{-1}.
\end{eqnarray*}
\end{lemma}
{\it Proof:} \cite[Corollary 2]{montvau}.
\section{Proof of Theorem \ref{th:3}}\setcounter{equation}{0}\setcounter{theorem}{0}\setcounter{lemma}{0}
 As in \cite[Proof of Theorem 3]{vau}, we square out the term $\Delta^{2}(n):$
\begin{eqnarray}
&&\sum\limits_{n\leq x\atop n\equiv N(mod\,v)}\Delta^{2}(n)= \sum\limits_{n\leq x\atop n\equiv N(mod\,v)}\Lambda^{2}(n)-
2\sum\limits_{r\leq R}\frac{\mu(r)}{\phi(r)}\sum\limits_{b=1\atop (b,r)=1}^{r}\sum\limits_{n\leq x\atop n\equiv N(mod\,v)}\Lambda(n)e\left(\frac{bn}{r}\right)+
\sum\limits_{n\leq x\atop n\equiv N(mod\,v)}|F_{R}(n)|^{2}.\nonumber\\
\label{eq:g1}
\end{eqnarray}
Using the prime number theorem in arithmetic progressions, we obtain
\begin{eqnarray}\label{eq:g2}
\sum\limits_{n\leq x\atop n\equiv N(mod\,v)}\Lambda^{2}(n)
=\delta(N,v)\frac{xL-x}{\phi(v)}+O\left(x\,exp\left(-cL^{-1/2}\right)\right).
\end{eqnarray}
 Applying the prime number theorem in arithmetic progressions once more, we see
\begin{eqnarray*}
&&\sum\limits_{b=1\atop (b,r)=1}^{r} \sum_{n\leq x\atop n\equiv N(mod\,v)} \Lambda(n)e(bn/r)
=\sum\limits_{n\leq x\atop {(n,r)=1\atop n\equiv N(mod\,v)}} \Lambda(n) \sum\limits_{b=1\atop (b,r)=1}^{r} e(bn/r) + O(RL^{2})\\
&=&\sum\limits_{n\leq x\atop {(n,r)=1\atop n\equiv N(mod\,v)}} \Lambda(n) \mu(r) + O(RL^{2})
 =\mu(r) \sum_{n\le x\atop n\equiv N(mod\,v)} \Lambda(n) + O(RL^{2}) \nonumber\\
&=& \delta(N,v)\frac{x\mu(r)}{\phi(v)} + O(x\,\exp(-cL^{1/2}))
,\end{eqnarray*} which implies that
\begin{eqnarray}\label{eq:g3a}
&&\sum\limits_{r\leq R}\frac{\mu(r)}{\phi(r)}\sum\limits_{b=1\atop (b,r)=1}^{r}\sum\limits_{n\leq x\atop n\equiv N(mod\,v)}\Lambda(n)e\left(\frac{bn}{r}\right)=\delta(N,v)\frac{x}{\phi(v)} \sum\limits_{r\leq R}\frac{\mu^{2}(r)}{\phi(r)}
+ O(x\,\exp(-cL^{1/2}/2)).\nonumber\\
\end{eqnarray} We know from \cite[Proof of Theorem 3]{vau} that
\begin{eqnarray}\label{eq:g3b}
\sum\limits_{r\leq R}\frac{\mu^{2}(r)}{\phi(r)}=\log\,R\,+\gamma+\sum\limits_{p\geq 2}\frac{\log\,p}{p(p-1)}+O\left(R^{-1/2}\right).
\end{eqnarray}
Inserting (\ref{eq:g3b}) into (\ref{eq:g3a}), we obtain
\begin{eqnarray}\label{eq:g3}
&&\sum\limits_{r\leq R}\frac{\mu(r)}{\phi(r)}\sum\limits_{b=1\atop (b,r)=1}^{r}\sum\limits_{n\leq x\atop n\equiv N(mod\,v)}\Lambda(n)e\left(\frac{bn}{r}\right)\nonumber\\
&=&\delta(N,v)\frac{x}{\phi(v)}\left(\log\,R\,+\gamma+\sum\limits_{p\geq 2 }\frac{\log\,p}{p(p-1)}\right)+O\left(\frac{x}{\phi(v)R^{1/2}}+x\,\exp(-cL^{1/2}/2)\right).\nonumber\\
\end{eqnarray}
We now estimate the third, most complicated term in (\ref{eq:g1}):
\begin{eqnarray}\label{eq:g4}
&&\sum\limits_{n\leq x\atop n\equiv N(mod\,v)}|F_{R}(n)|^{2}=\sum\limits_{r\leq R}\frac{\mu(r)}{\phi(r)}
\sum\limits_{b=1\atop (b,r)=1}^{r}\sum\limits_{r_{1}\leq R}\frac{\mu(r_{1})}{\phi(r_{1})}\sum\limits_{b_{1}=1\atop (b_{1},r_{1})=1}^{r_{1}}
\frac{\mu(r_{1})}{\phi(r_{1})}
\sum\limits_{n\leq x \atop n\equiv N(mod\,v)}e\left(\frac{bn}{r}-\frac{b_{1}n}{r_{1}}\right)\nonumber\\
&=:& A_{N}+2B_{N}+C_{N},
\end{eqnarray} where
\begin{eqnarray}\label{eq:g5a}
A_{N}&=& \left(\frac{x}{v}+O(1)\right)\sum\limits_{r|v\atop r\leq R}\frac{\mu(r)}{\phi(r)}
\sum\limits_{b=1\atop (b,r)=1}^{r}\sum\limits_{r_{1}|v\atop r_{1}\leq R}\frac{\mu(r_{1})}{\phi(r_{1})}\sum\limits_{b_{1}=1\atop (b_{1},r_{1})=1}^{r_{1}}
e\left(\frac{bN}{r}-\frac{b_{1}N}{r_{1}}\right)\nonumber\\
&=&\left(\frac{x}{v}+O(1)\right)\left(\sum\limits_{r|v\atop r\leq R}\frac{\mu(r)C_{r}(N)}{\phi(r)}\right)^{2},\\
\label{eq:g5}B_{N}&=&
\sum\limits_{r| v\atop r\leq R}\frac{\mu(r)}{\phi(r)}
\sum\limits_{b=1\atop (b,r)=1}^{r}\sum\limits_{r_{1}\nmid v\atop r_{1}\leq R}\frac{\mu(r_{1})}{\phi(r_{1})}\sum\limits_{b_{1}=1\atop (b_{1},r_{1})=1}^{r_{1}}
\sum\limits_{n\leq x\atop n\equiv N(mod\,v)}e\left(\frac{bn}{r}-\frac{b_{1}n}{r_{1}}\right)\nonumber\\
&=&\sum\limits_{r| v\atop r\leq R}\frac{\mu(r)}{\phi(r)}C_{r}(N)
\rho^{*}(x,v,N)\ll \tau(v)|\rho^{*}(x,v,N)|,\end{eqnarray}
where $C_{r}(N)$ denotes the Ramanujuan sum and
\begin{eqnarray}\label{eq:rho2}
\rho^{*}(x,v,N)=\sum\limits_{r_{1}\nmid v\atop r_{1}\leq R}
\frac{\mu(r_{1})}{\phi(r_{1})}\sum\limits_{b_{1}=1\atop (b_{1},r_{1})=1}^{r_{1}}
\sum\limits_{n\leq x\atop n\equiv N(mod\,v)}e\left(-\frac{b_{1}n}{r_{1}}\right).
\end{eqnarray}
\begin{eqnarray}\label{eq:g5C} C_{N}&=&
\sum\limits_{r\nmid v\atop r\leq R}\frac{\mu(r)}{\phi(r)}
\sum\limits_{b=1\atop (b,r)=1}^{r}\sum\limits_{r_{1}\nmid v\atop r_{1}\leq R}\frac{\mu(r_{1})}{\phi(r_{1})}\sum\limits_{b_{1}=1\atop (b_{1},r_{1})=1}^{r_{1}}
\sum\limits_{n\leq x\atop n\equiv N(mod\,v)}e\left(\frac{bn}{r}-\frac{b_{1}n}{r_{1}}\right)\nonumber\\
&=&C_{N,1}+C_{N,2},\end{eqnarray}
where
\begin{eqnarray}
\label{eq:C1}C_{N,1}&=&
\left(\frac{x}{v}+O(1)\right)\sum\limits_{r\nmid v\atop r\leq R}\frac{\mu^{2}(r)}{\phi(r)}
=\left(\frac{x}{v}+O(1)\right)\left(\sum\limits_{ r\leq R}-\sum\limits_{r| v}+\sum\limits_{r|v \atop r> R}\right)\frac{\mu^{2}(r)}{\phi(r)}\nonumber\\
&=&\left(\frac{x}{v}+O(1)\right)\left(\sum\limits_{r\leq R}\frac{\mu^{2}(r)}{\phi(r)}-\frac{v}{\phi(v)}+O\left(\sum\limits_{r| v\atop r>R}\frac{1}{\phi(r)}\right)\right),\nonumber\\
&=&\frac{x}{v}\sum\limits_{ r\leq R}\frac{\mu^{2}(r)}{\phi(r)}-\frac{x}{\phi(v)}+O\left(\frac{x\tau(v)\log\,R}{vR^{}}\right),\\
\label{eq:C2}C_{N,2}&=&\sum\limits_{r\nmid v\atop r\leq R}\frac{\mu(r)}{\phi(r)}
\sum\limits_{b=1\atop (b,r)=1}^{r}e\left(\frac{Nb}{r}\right)\sum\limits_{r_{1}\nmid v\atop r_{1}\leq R}\frac{\mu(r_{1})}{\phi(r_{1})}\sum\limits_{b_{1}=1\atop{ (b_{1},r_{1})=1\atop (r,b)\ne (r_{1},b_{1})}}^{r_{1}}e\left(\frac{-Nb_{1}}{r_{1}}\right)
\sum\limits_{s\leq (x-N)/v}e\left(\frac{sv(br_{1}-b_{1}r)}{rr_{1}}\right)\nonumber\\
&=&\sum\limits_{r\nmid v\atop r\leq R}\frac{\mu(r)}{\phi(r)}
\sum\limits_{b=1\atop (b,r)=1}^{r}e\left(\frac{(N+v)b}{r}\right)\sum\limits_{r_{1}\nmid v\atop r_{1}\leq R}\frac{\mu(r_{1})}{\phi(r_{1})}\sum\limits_{b_{1}=1\atop{ (b_{1},r_{1})=1\atop (r,b)\ne (r_{1},b_{1})}}^{r_{1}}e\left(\frac{-(N+v)b_{1}}{r_{1}}\right)\nonumber\\
&\times &\frac{ e\left(\lfloor (x-N)/v \rfloor\left(\frac{vb}{r}-\frac{vb_{1}}{r_{1}}\right)\right)-1}
{e\left(\frac{vb}{r}-\frac{vb_{1}}{r_{1}}\right)-1}\nonumber\\
&=&\sum\limits_{r\nmid v\atop r\leq R}\frac{\mu(r)}{\phi(r)}
\sum\limits_{b=1\atop (b,r)=1}^{r}e\left(\frac{(N+v/2)b}{r}\right)
e\left(\lfloor (x-N)/v \rfloor\frac{vb}{2r}\right)\nonumber\\
&\times &
\sum\limits_{r_{1}\nmid v\atop r_{1}\leq R}\frac{\mu(r_{1})}{\phi(r_{1})}\sum\limits_{b_{1}=1\atop{ (b_{1},r_{1})=1\atop (r,b)\ne (r_{1},b_{1})}}^{r_{1}}e\left(\frac{-(N+v/2)b_{1}}{r_{1}}\right) e\left(-\lfloor (x-N)/v \rfloor\frac{vb}{2r_{1}}\right)\nonumber\\
&\times &\frac{ \sin\left(\pi \lfloor (x-N)/v \rfloor\left(\frac{vb}{r}-\frac{vb_{1}}{r_{1}}\right)\right)}
{\sin\left(\pi\left(\frac{vb}{r}-\frac{vb_{1}}{r_{1}}\right)\right)}.\nonumber
\end{eqnarray}  Using Lemma \ref{le:hilbert} with $v_{r}:=v_{r,b}=\frac{\mu(r)}{\phi(r)}e\left(\frac{(N+v/2)b}{r}\right)
e\left(\lfloor (x-N)/v \rfloor\frac{vb}{2r}\right),$ $x_{r}:=x_{r,b}=\frac{vb}{r},$ and $\delta\gg R^{-2},$ we obtain
\begin{eqnarray}\label{eq:C2a}
C_{N,2}\ll R^{2}\log\,R.
\end{eqnarray}
Combining (\ref{eq:g5C}), (\ref{eq:C1}) and (\ref{eq:C2a}), we see
\begin{eqnarray}\label{eq:C3}
C_{N}&=&\frac{x}{v}\sum\limits_{r\leq R}\frac{\mu^{2}(r)}{\phi(r)}-\frac{x}{\phi(v)}+O\left(\frac{x\tau(v)\log\,R}{vR}+R^{2}\log\,R\right).
\end{eqnarray}
To evaluate $A_{N},$ we see that for squarefree $v$ and $N\geq 1,$
\begin{eqnarray}\label{eq:g11}
&&\sum\limits_{r|v}\frac{\mu(r)C_{r}(N)}{\phi(r)}=\sum\limits_{r|v}\frac{\mu(r)}{\phi(r)}\sum\limits_{d|(N,r)}
d\mu\left(\frac{r}{d}\right)
=\sum\limits_{d|(N,v)}d\sum\limits_{d|r\atop r|v}\frac{\mu(r)}{\phi(r)}\mu\left(\frac{r}{d}\right)\nonumber\\
&=&\sum\limits_{d|(N,v)}d\sum\limits_{ h|\frac{v}{d}}\frac{\mu(hd)}{\phi(hd)}\mu\left(h\right)=
\sum\limits_{d|(N,v)}\frac{\mu(d)d}{\phi(d)}\sum\limits_{ h|\frac{v}{d}}\frac{\mu^{2}(h)}{\phi(h)}=
\sum\limits_{d|(N,v)}\frac{\mu(d)d}{\phi(d)\phi\left(\frac{v}{d}\right)}\frac{v}{d}\nonumber\\
&=&\frac{v}{\phi(v)}\sum\limits_{d|(N,v)}\mu(d)=\frac{v}{\phi(v)}\delta(N,v).
\end{eqnarray}
For $N=0,$ we see similarly,
\begin{eqnarray}\label{eq:g11aa}
&&\sum\limits_{r|v}\frac{\mu(r)C_{r}(N)}{\phi(r)}=\sum\limits_{r|v}\mu(r)=\frac{v}{\phi(v)}\delta(N,v).
\end{eqnarray}
Further, we note that
\begin{eqnarray}\label{eq:g12}
&&\left|\sum\limits_{r|v\atop r>R}\frac{\mu(r)C_{r}(N)}{\phi(r)}\right|\leq \sum\limits_{r|v\atop r>R}1
\leq \tau(v).
\end{eqnarray}
Combining (\ref{eq:g5a}) and (\ref{eq:g11}) - (\ref{eq:g12}), we see
\begin{eqnarray}\label{eq:g13}
&&A_{N}=\left(\frac{x}{v}+O(1)\right)\left(\frac{v}{\phi(v)}\delta(N,v)+O\left(\left|\sum\limits_{r|v\atop r>R}\frac{\mu(r)C_{r}(N)}{\phi(r)}\right|\right)\right)^{2}\nonumber\\
&=&x\frac{v}{\phi^{2}(v)}\delta(N,v)+O\left(\log\,v+\frac{x(\log\,v+\tau(v))}{v}\sum\limits_{r|v\atop r>R}1\right)\nonumber\\
&=& x\frac{v}{\phi^{2}(v)}\delta(N,v)+O\left(\frac{x(\log\,v+\tau(v))}{v}\sum\limits_{r|v\atop r>R}1\right).
\end{eqnarray}
 We know from \cite[Proof of Theorem 1]{vau} that for $\rho^{*}(x,v,N)$ as defined in (\ref{eq:rho2}):
\begin{eqnarray}\label{eq:g6}
\rho^{*}(x,v,N)\ll R.
\end{eqnarray}
Combining (\ref{eq:g5}) and (\ref{eq:g6}), we get
\begin{eqnarray}\label{eq:g14}
B_{N}\ll \tau(v)R.
\end{eqnarray} Combining (\ref{eq:g4}), (\ref{eq:C3}), (\ref{eq:g13}), and (\ref{eq:g14}), and subsequently applying (\ref{eq:g3b}), we obtain
\begin{eqnarray}\label{eq:g15}
&&\sum\limits_{n\leq x\atop n\equiv N(mod\,v)}|F_{R}(n)|^{2}=x\frac{v}{\phi^{2}(v)}\delta(N,v)+\frac{x}{v}\sum\limits_{r\leq R}\frac{\mu^{2}(r)}{\phi(r)}-\frac{x}{\phi(v)}\nonumber\\
&+&O\left(\frac{x\tau(v)\log \,R}{vR}+R^{2}\log\,R+\tau(v)R+\frac{x(\log\,v+\tau(v))}{v}\sum\limits_{r|v\atop r>R}1\right)\nonumber\\
&=&x\frac{v}{\phi^{2}(v)}\delta(N,v)+\frac{x}{v}\left(\log\,R\,+\gamma+\sum\limits_{p\geq }\frac{\log\,p}{p(p-1)}\right)-\frac{x}{\phi(v)}\nonumber\\
&+&O\left(\frac{x\tau(v)}{vR^{1/2}}+R^{2}\log\,R+\tau(v)R+\frac{x(\log\,v+\tau(v))}{v}\sum\limits_{r|v\atop r>R}1\right).
\end{eqnarray}
Finally, we derive Theorem \ref{th:3} from (\ref{eq:g1}), (\ref{eq:g2}), (\ref{eq:g3}), and (\ref{eq:g15}).
\section{Proof of Theorem \ref{th:4}}\setcounter{equation}{0}\setcounter{theorem}{0}\setcounter{lemma}{0}
For $Q_{1}=x/R,$ we see by arguing as in \cite[Proof of Theorem 4]{vau} that
\begin{eqnarray}\label{eq:RRR}
\sum\limits_{d\leq Q_{1}}\sum\limits_{b=1\atop (N-b,d)=1}^{d}\left(\vartheta(x,d,b)-\rho(x,d,b)\right)^{2} \ll x^{2}R^{-1}L^{2}.
\end{eqnarray}
Thus, for $Q_{1}<Q$ we only need to consider the expression
 \begin{eqnarray}\label{eq:th41}
 &&\sum\limits_{Q_{1}\leq d\leq Q}\sum\limits_{b=1\atop (N-b,d)=1}^{d}\left(\vartheta(x,d,b)-\rho(x,d,b)\right)^{2}\nonumber\\
 &=&
\sum\limits_{Q_{1}<d<Q}\sum\limits_{v|d}\mu(v)\sum\limits_{b=1\atop v|N-b}^{d}\left(\vartheta(x,d,b)-\rho(x,d,b)\right)^{2}\nonumber\\
&=&\sum\limits_{v<Q}\mu(v)\sum\limits_{Q_{1}<d<Q\atop v|d}\sum\limits_{b=1\atop v|N-b}^{d}\left(\vartheta(x,d,b)-\rho(x,d,b)\right)^{2}.
\end{eqnarray}
We separately treat the cases $v\leq R^{5}$ and $v>R^{5}.$ In the second case we see by using a trivial estimate for $\vartheta(x,q,a)$ and applying Lemma \ref{le:hy} to estimate $\rho(x,q,a),$
\begin{eqnarray}\label{eq:auv2}
&&\sum\limits_{v>R^{5}}\mu(v)\sum\limits_{Q_{1}<d<Q\atop v|d}\sum\limits_{b=1\atop v|N-b}^{d}\left(\vartheta(x,d,b)-\rho(x,d,b)\right)^{2}\nonumber\\
&\ll &
\sum\limits_{v>R^{5}}\sum\limits_{Q_{1}<d<Q\atop v|d}\sum\limits_{b=1\atop v|N-b}^{d}\vartheta^{2}(x,d,b)+\rho^{2}(x,d,b)\nonumber\\
&\ll &x^{2}L^{4}\sum\limits_{v>R^{5}}\sum\limits_{Q_{1}<d<Q\atop v|d}\frac{\tau(d)}{d^{2}}\sum\limits_{b=1\atop v|N-b}^{d}1+R^{2}\sum\limits_{v>R^{5}}\sum\limits_{Q_{1}<d<Q\atop v|d}\sum\limits_{b=1\atop v|N-b}^{d}1\nonumber\\
&\ll & x^{2}L^{4}\sum\limits_{v>R^{5}}\frac{1}{v}\sum\limits_{Q_{1}<d<Q\atop v|d}\frac{\tau(d)}{d}
+R^{2}\sum\limits_{v>R^{5}}\frac{1}{v}\sum\limits_{Q_{1}<d<Q\atop v|d}d\nonumber\\
&\ll&  x^{2}L^{6}\sum\limits_{v>R^{5}}\frac{1}{v^{2}}+R^{2}Q^{2}\sum\limits_{v>R^{5}}\frac{1}{v^{2}}\ll   x^{2}R^{-3}.
\end{eqnarray}
In the case $v\leq R^{5},$ we see for a fixed $v,$
\begin{eqnarray}\label{eq:uv4}
&&\sum\limits_{Q_{1}<d<Q\atop v|d}\sum\limits_{b=1\atop v|N-b}^{d}\left(\vartheta(x,d,b)-\rho(x,d,b)\right)^{2}=\sum\limits_{\frac{Q_{1}}{v}<m<\frac{Q}{v}}\sum\limits_{b=1\atop v|N-b}^{mv}\left(\vartheta(x,mv,b)-\rho(x,mv,b)\right)^{2}
\nonumber\\
&=&\sum\limits_{\frac{Q_{1}}{v}<m<\frac{Q}{v}}\sum\limits_{n_{1},n_{2}\leq x\atop{ n_{1}\equiv n_{2}(mod\,mv)\atop n_{1}\equiv n_{2}\equiv N(mod\,v)}}\Delta(n_{1})\Delta(n_{2}),\nonumber\\
&=&\sum\limits_{\frac{Q_{1}}{v}<m<\frac{Q}{v}}\sum\limits_{n\leq x\atop{ n\equiv N(mod\,v)}}\Delta^{2}(n)+2\sum\limits_{\frac{Q_{1}}{v}<m<\frac{Q}{v}}\sum\limits_{n_{1}<n_{2}\leq x\atop{ n_{1}\equiv n_{2}(mod\,mv)\atop n_{1}\equiv n_{2}\equiv N(mod\,v)}}\Delta(n_{1})\Delta(n_{2})\nonumber\\
&:=&E_{v}+2F_{v}.
\end{eqnarray}
Applying Theorem \ref{th:3}, we see
\begin{eqnarray}\label{eq:EE}
E_{v}&=&\sum\limits_{\frac{Q_{1}}{v}<m<\frac{Q}{v}}\left(\delta(N,v)\frac{x}{\phi(v)}
\left(\log\,x-2\log\,R-c_{1}\right)+\frac{x}{v}\left(\log\,R+c_{2}\right)+\delta(N,v)x\frac{v}{\phi^{2}(v)}-\frac{x}{\phi(v)}\right)\nonumber\\
&+&O\left(Q\left(\frac{x\,exp\left(-cL^{1/2}\right)}{v}+\frac{x\tau(v)}{v^{2}R^{1/2}}+\frac{x}{v\phi(v)R^{1/2}}+\frac{R^{2}\log\,R}{v}
+\frac{\tau(v)R}{v}+\frac{x(\log\,v+\tau(v))}{v^{2}}\sum\limits_{r|v\atop r>R}1\right)\right)\nonumber\\
&=&\delta(N,v)\frac{xQ}{v\phi(v)}
\left(\log\,x-2\log\,R-c_{1}\right)+\frac{xQ}{v^{2}}\left(\log\,R+c_{2}\right)+\delta(N,v)\frac{xQ}{\phi^{2}(v)}-\frac{xQ}{v\phi(v)}\nonumber\\
&+&O\left(\frac{xQ}{v^{3/2}R^{1/2}}+\frac{xQ(\log\,v+\tau(v))}{v^{2}}\sum\limits_{r|v\atop r>R}1\right),
\end{eqnarray} from which we derive
\begin{eqnarray}\label{eq:g20}
&&\sum\limits_{v\leq R^{5}}\mu(v)E_{v}\nonumber\\
&=&\sum\limits_{v\leq R^{5}\atop (v,N)=1}\mu(v)\left(\frac{xQ}{v\phi(v)}\left(\log\,x-2\log\,R-c_{1}\right)+\frac{xQ}{\phi^{2}(v)}\right)+
\sum\limits_{v\leq R^{5}}\mu(v)\left(\frac{xQ}{v^{2}}\left(\log\,R+c_{2}\right)-\frac{xQ}{v\phi(v)}\right)\nonumber\\
&+&O\left(\sum\limits_{v\leq R^{5}}\frac{xQ}{v^{3/2}R^{1/2}}+\sum\limits_{v\leq R^{5}}\frac{xQ(\log\,v+\tau(v))}{v^{2}}\sum\limits_{r|v\atop r>R}1\right)\nonumber\\
&=&\sum\limits_{v\geq 1\atop (v,N)=1}\mu(v)\left(\frac{xQ}{v\phi(v)}\left(\log\,x-2\log\,R-c_{1}\right)+\frac{xQ}{\phi^{2}(v)}\right)+
\sum\limits_{v\geq 1}\mu(v)\left(\frac{xQ}{v^{2}}\left(\log\,R+c_{2}\right)-\frac{xQ}{v\phi(v)}\right)\nonumber\\
&+&O\left(\frac{xQ}{R^{1/2}}+
xQ\sum\limits_{R<r\leq R^{5}}\sum\limits_{v\leq R^{5}\atop r|v}\frac{1}{v^{2-\epsilon}}\right)\nonumber\\
&=&xQ\prod\limits_{p\geq 2\atop (p,N)=1}\left(1-\frac{1}{p(p-1)}\right)\left(\log\,x-2\log\,R-c_{1}\right)+
xQ\prod\limits_{p\geq 2\atop (p,N)=1}\left(1-\frac{1}{(p-1)^{2}}\right)\nonumber\\
&+&xQ\prod\limits_{p\geq 2}\left(1-\frac{1}{p^{2}}\right)\left(\log\,R+c_{2}\right)-xQ\prod\limits_{p\geq 2}\left(1-\frac{1}{p(p-1)}\right)+O\left(\frac{xQ}{R^{1/2}}+
xQ\sum\limits_{R<r\leq R^{5}}\frac{1}{r^{2-\epsilon}}\right)\nonumber\\
&=&xQ\prod\limits_{p\geq 2\atop (p,N)=1}\left(1-\frac{1}{p(p-1)}\right)\left(\log\,x-(\log\,R){t(N)}\right)\nonumber\\
&+&xQ\Biggl(-\prod\limits_{p\geq 2\atop (p,N)=1}\left(1-\frac{1}{p(p-1)}\right)c_{1}+
\prod\limits_{p\geq 2\atop (p,N)=1}\left(1-\frac{1}{(p-1)^{2}}\right)+\prod\limits_{p\geq 2}\left(1-\frac{1}{p^{2}}\right)c_{2}\nonumber\\
&&-\prod\limits_{p\geq 2}\left(1-\frac{1}{p(p-1)}\right)\Biggr)+O\left(\frac{xQ}{R^{1/2}}\right),
\end{eqnarray}
where
\begin{eqnarray*}
t(N) &=&
\left(2-\prod\limits_{p\geq 2\atop (p,N)=1}\left(1-\frac{1}{p(p-1)}\right)^{-1}\prod\limits_{p\geq 2}\left(1-\frac{1}{p^{2}}\right)\right)\geq 1.
\end{eqnarray*}
For $F_{v},$ we see
\begin{eqnarray}\label{eq:th45}
F_{v}&=&B_{v,\frac{Q_{1}}{v}}-B_{v,\frac{Q}{v}},
\end{eqnarray}
where
\begin{eqnarray*}
B_{v,u}=\sum\limits_{u<m\leq x}\sum\limits_{n_{1}<n_{2}\leq x\atop {n_{2}-n_{1}\equiv N(mod\,vm)\atop n_{1}\equiv n_{2}\equiv N(mod\,v)}}\Delta(n_{1})\Delta(n_{2}).
\end{eqnarray*}
To calculate $F_{v,u},$ we first note that the summation condition $m\leq x$ is superfluous as $n_{1},n_{2}\leq x.$ Using this fact and applying the divisor switching trick from \cite{hooley1}, we see
\begin{eqnarray} \label{eq:uv7}
B_{v,u}&=&\sum\limits_{u<m}\sum\limits_{n_{1}\leq x\atop  n_{1}\equiv N(mod\,v)}\sum\limits_{n_{1}<n_{2}\leq x\atop {n_{2}\equiv n_{1}(mod\,vm)\atop  n_{2}\equiv N(mod\,v)}}\Delta(n_{1})\,\Delta(n_{2})\nonumber\\
&=&\sum\limits_{l\leq x/vu}\sum\limits_{u<m}\sum\limits_{n_{1}\leq x\atop  n_{1}\equiv N(mod\,v)}\sum\limits_{lvu+n_{1}<n_{2}\leq x\atop {n_{2}-n_{1}=lvm\atop  n_{2}\equiv N(mod\,v)}}\Delta(n_{1})\,\Delta(n_{2})\nonumber\\
&=&\sum\limits_{l\leq x/vu}\sum\limits_{n_{1}\leq x-lvu\atop  n_{1}\equiv N(mod\,v)}\sum\limits_{lvu+n_{1}<n_{2}\leq x\atop {n_{2}\equiv n_{1}(mod \,lv)\atop  n_{2}\equiv N(mod\,v)}}\Delta(n_{1})\,\Delta(n_{2})\nonumber\\
&=&\sum\limits_{l\leq x/vu}\sum\limits_{a=1\atop{  a\equiv N(mod\,v)}}^{lv}\sum\limits_{n_{1}<x-lvu\atop n_{1}\equiv a(mod\,lv)}\Delta(n_{1})\sum\limits_{lvu+n_{1}<n_{2}\leq x\atop n_{2}\equiv a(mod\,lv)}\Delta(n_{2}).
\end{eqnarray}
Applying Lemma \ref{le:th2vau} to the sum over $n_{2},$ we see
\begin{eqnarray}
\sum\limits_{lvu+n_{1}<n_{2}\leq x\atop n_{2}\equiv a(mod\,lv)}\Delta(n_{2}) \ll x\,exp\left(-cL^{1/2}\right).
\end{eqnarray}
Applying Theorem \ref{th:3} and the Cauchy inequality, we see
\begin{eqnarray}\label{eq:h6}
\sum\limits_{n_{1}<x-lvu\atop n_{1}\equiv a(mod\,lv)}|\Delta(n_{1})|\ll xL.
\end{eqnarray}
Combining (\ref{eq:uv7}) - (\ref{eq:h6}), we obtain
\begin{eqnarray*}\label{eq:u67}
&&B_{v,u}\ll  x^{2}L\,exp\left(-cL^{1/2}\right)\sum\limits_{l\leq x/vu}\sum\limits_{a=1\atop{  a\equiv N(mod\,v)}}^{lv}  1
\ll  x^{2}L\,exp\left(-cL^{1/2}/2\right),
\end{eqnarray*} which together with (\ref{eq:th45}) implies
\begin{eqnarray}\label{eq:u677}
\sum\limits_{v\leq R^{5}}\mu(v)F_{v}\ll x^{2}\,exp\left(-cL^{1/2}/4\right).
\end{eqnarray}
Combining (\ref{eq:RRR}) - (\ref{eq:g20}) and (\ref{eq:u677}), we derive Theorem \ref{th:4}.

\section{Proof of Theorem \ref{th:5}}\setcounter{equation}{0}\setcounter{theorem}{0}\setcounter{lemma}{0}
Arguing as in (\ref{eq:RRR}) and (\ref{eq:th41}) in the proof of Theorem \ref{th:4}, we can limit ourselves to the expression
 \begin{eqnarray}\label{eq:hth41}
 &&\sum\limits_{Q_{1}\leq d\leq Q}\sum\limits_{b=1\atop (b,d)=1}^{d}\left(\vartheta(x,d,b)-\rho(x,d,b)\right)^{2}\nonumber\\
&=&\sum\limits_{v<Q}\mu(v)\sum\limits_{Q_{1}<d<Q\atop v|d}\sum\limits_{b=1\atop v|b}^{d}\left(\vartheta(x,d,b)-\rho(x,d,b)\right)^{2}.
\end{eqnarray}
We separately treat the cases $v\leq R^{5}$ and $v>R^{5}.$ In the second case, we see by arguing as in (\ref{eq:auv2}),
\begin{eqnarray}\label{eq:h32}
&&\sum\limits_{v>R^{5}}\mu(v)\sum\limits_{Q_{1}<d<Q\atop v|d}\sum\limits_{b=1\atop v|b}^{d}\left(\vartheta(x,d,b)-\rho(x,d,b)\right)^{2}
\ll  x^{2}R^{-3}.
\end{eqnarray}
In the case $v\leq R^{5},$ we see by arguing as in (\ref{eq:uv4}) for a fixed $v,$
\begin{eqnarray}\label{eq:uv4hh}
&&\sum\limits_{Q_{1}<d<Q\atop v|d}\sum\limits_{b=1\atop v|b}^{d}\left(\vartheta(x,d,b)-\rho(x,d,b)\right)^{2}=\sum\limits_{\frac{Q_{1}}{v}<m<\frac{Q}{v}}\sum\limits_{b=1\atop v|b}^{mv}\left(\vartheta(x,mv,b)-\rho(x,mv,b)\right)^{2}
\nonumber\\
&=&\sum\limits_{\frac{Q_{1}}{v}<m<\frac{Q}{v}}\sum\limits_{n\leq x\atop{ n\equiv 0(mod\,v)}}\Delta^{2}(n)+2\sum\limits_{\frac{Q_{1}}{v}<m<\frac{Q}{v}}\sum\limits_{n_{1}<n_{2}\leq x\atop {n_{2}\equiv n_{1}(mod\,mv)\atop n_{1}\equiv n_{2}\equiv 0(mod\,v)}}\Delta(n_{1})\Delta(n_{2})\nonumber\\
&:=&E_{v}+2H_{v},
\end{eqnarray} where $E_{v}$ is as given in (\ref{eq:EE}) with $N\equiv 0.$
Following the argument in (\ref{eq:g20}) and recalling the definition of $\delta(N,v)$ in Theorem \ref{th:3}, we obtain
\begin{eqnarray}\label{eq:h20}
&&\sum\limits_{v\leq R^{5}}\mu(v)E_{v}\nonumber\\
&=&xQ\left(\log\,x-2\log\,R-c_{1}+1\right)
+\sum\limits_{v\geq 1}\left(
\frac{xQ}{v^{2}}\left(\log\,R+c_{2}\right)-\frac{xQ}{v\phi(v)}\right)+O\left(\frac{xQ}{R^{1/2}}\right)\nonumber\\
&=&xQ\left(\log\,x-2\log\,R-c_{1}+1\right)
+
xQ\prod\limits_{p\geq 2}\left(1-\frac{1}{p^{2}}\right)\left(\log\,R+c_{2}\right)-xQ\prod\limits_{p\geq 2}\left(1-\frac{1}{p(p-1)}\right)\nonumber\\
&+&O\left(\frac{xQ}{R^{1/2}}\right)\nonumber\\
&=&xQ\left(\log\,x-(\log\,R)\left(2-\zeta^{-1}(2)\right)\right)
+xQ\left(-c_{1}+1
-\prod\limits_{p\geq 2}\left(1-\frac{1}{p(p-1)}\right)+\prod\limits_{p\geq 2}\left(1-\frac{1}{p^{2}}\right)c_{2}\right)\nonumber\\
&+&O\left(\frac{xQ}{R^{1/2}}\right),
\end{eqnarray}
For $H_{v},$ we see by arguing as in (\ref{eq:th45}) - (\ref{eq:u677}),
\begin{eqnarray}\label{eq:hu677}
\sum\limits_{v\leq R^{5}}\mu(v)F_{v}\ll x^{2}\,exp\left(-cL^{1/2}/4\right).
\end{eqnarray} Combining (\ref{eq:hth41}) -  (\ref{eq:hu677}), we derive Theorem \ref{th:5}.

\end{document}